\documentclass[preprint,review,12pt]{elsarticle}

\usepackage{graphics}
\usepackage{graphicx}
\usepackage{epsfig}

\usepackage{amssymb}

\begin{document}

\begin{frontmatter}



\title{Common fixed
points and endpoints of multi-valued  generalized weak contraction
mappings}


\author{Congdian Cheng}

\address{College of Mathematics and Systems
Science, Shenyang Normal University,
  Shenyang,  110034,   China}

\begin{abstract}
Let $(X, d)$ be a complete metric space, and let $S, T :
X\rightarrow
  CB(X)$ be a duality of multi-valued generalized
weak contraction mappings or a duality of generalized $\varphi$-weak
contraction mappings. We discuss the common fixed points and
endpoints of the two kinds of multi-valued   weak mappings.  Our
results extend and improve some results given by  Daffer and Kaneko
(1995), Rouhani and Moradi (2010), and Moradi and Khojasteh (2011).
\end{abstract}


\end{frontmatter}


\section{Introduction}
\label{1}Let $(X, d)$ be a metric space and $CB(X)$ denote the class
of closed and bounded subsets of X.  Also let $S, T : X\rightarrow
2^X$ be a multi-valued mapping. A point x is called a fixed point of
$T$ if $x\in Tx$. Define $Fix(T) = \{x\in X : x \in Tx\}$. An
element $x \in X$ is said to be an endpoint (or stationary point) of
a multi-valued mapping $T$ if $Tx = \{x\}$. We denote the set of all
endpoints of $T$ by $End(T)$.

A bivariate mapping $\phi: X\times X\rightarrow[0, +\infty)$ is
called compactly positive if $\inf\{\phi(x,y): a\leq d(x,y)\leq
b\}>0$ for each finite interval $[a, b]\subseteq (0, +\infty)$. A
mapping $T : X \rightarrow CB(X)$ is called weakly contractive if
there exists a compactly positive mapping $\phi$
 such that
$$H(Tx, Ty)
\leq d(x, y)-\phi(x, y)$$ for each $x, y \in X$, where
$$H(A, B) :=\max\{ \sup\limits_{x\in B}
d(x, A), \sup\limits_{x\in A} d(x, B)\},$$ denoting the Hausdorff
metric on $CB(X)$. (see [1].)

A mapping $T : X \rightarrow CB(X)$  is called an generalized
$\varphi$-weak contraction if there exists a map $\varphi: [0,
+\infty)\rightarrow[0, +\infty)$ with $\varphi(0)=0$ and
$\varphi(t)<t$ for all $t>0$ such that
$$H(Sx, Ty)
\leq \varphi(N(x, y))$$
 for all $x, y \in X$, where
$$N(x, y)
:= \max \{d(x, y), d(x, Tx), d(y, Ty), \frac{d(x, Ty) + d(y,
Tx)}{2}\}. $$

Two mappings $S, T : X \rightarrow CB(X)$ ($S, T : X \rightarrow X$)
are called a duality of generalized weak contractions if there
exists a bivariate mapping $\alpha: X\times X\rightarrow[0, 1)$ such
that
$$H(Sx, Ty)
\leq \alpha(x, y)M(x, y)$$ for all $x, y \in X$ (or equivalently, if
there exists a bivariate mapping $\phi: [0, +\infty)\rightarrow[0,
+\infty)$ with $\phi(0)=0$ and $\phi(t)>0$ for all $t>0$
 such that
$$H(Sx, Ty)
\leq M(x, y)-\phi(x, y)$$ for each $x, y \in X$), where
$$\begin{array}{rcl}M(x, y): =
\max\{d(x,y), d(x, Sx), d(y, Ty), \frac{d(x, Ty)+d(y,
Sx)}{2}\}.\end{array}$$

Also, two mappings $S, T : X \rightarrow CB(X)$ are called a duality
of generalized $\varphi$-weak contractions if there exists a
bivariate mapping $\varphi: [0, +\infty)\rightarrow[0, +\infty)$
with $\varphi(0)=0$ and $\varphi(t)<t$ for all $t>0$ such that
$$H(Sx, Ty)
\leq \varphi(M(x, y))$$
 for all $x, y \in X$ (or
equivalently, if there exists a bivariate mapping $\varphi: [0,
+\infty)\rightarrow[0, +\infty)$ with $\varphi(0)=0$ and
$\varphi(t)>0$ for all $t>0$ such that
$$H(Sx, Ty)
\leq M(x, y)-\varphi(M(x, y))$$ for all $x, y \in X$).

A mapping $T : X \rightarrow CB(X)$ has the approximate endpoint
property if $$\inf\limits_{x\in X}\sup\limits_{y\in Tx}d(x, y) =
0.$$

The fixed points   for multi-valued  contraction mappings have been
the subject of the research area on fixed points for more than forty
years, for example, see [1-5] and the references therein. The
investigation of endpoints of multi-valued mappings  was made as
early as 30 years ago, and has received great attention in recent
years, see e.g. [5-10]. Among other studies, several important
results  related closely to the present work are as follows.

First, in the following theorem, Nadler [2] (1969) extended the
Banach contraction principle to multi-valued mappings.\\

\noindent{\bf Theorem 1.1.} \textit{Let $(X, d)$ be a complete
metric space. Suppose that $T : X \rightarrow CB(X)$ is a
contraction mapping in the sense that for some $0\leq \alpha<1,
H(Tx, Ty) \leq \alpha d(x, y)$ for all $x, y \in X$. Then there
exists a point $x \in X$ such that $x \in Tx$.}\\

\noindent Then, Daffer and Kaneko [1] (1995) proved the next Theorem
1.2 and
Theorem 1.3.\\

\noindent{\bf Theorem 1.2} ([1, Theorem 3.3]) \textit{Let $(X, d)$
be a complete metric space. Suppose that $T : X \rightarrow CB(X)$
be such that $H(Tx, Ty) \leq \alpha N(x, y)$ for $0\leq\alpha<1$,
for all $x, y \in X$. If $x\rightarrow d(x, Tx)$ is lower
semicontinuous (l.s.c.), then there
exists a point $x_0\in X$ such that $x_0\in Tx_0$.}\\

\noindent{\bf Theorem 1.3} ([1, Theorem 2.3]). \textit{Let $(X, d)$
be a complete metric and $T : X \rightarrow CB(X)$ weakly
contractive. Assume that
$$\liminf
\limits_{ \beta\rightarrow 0}\frac{\lambda(\alpha, \beta)}{\beta} >
0 \hspace{3mm}(0<\alpha\leq\beta),$$ where $\lambda(\alpha,
\beta)=\inf\{\phi(x, y)|x, y\in X, \alpha\leq d(x, y)\leq\beta\}$
for each finite interval $[\alpha, \beta]\subset(0, \infty)$. Then
$T$
has a fixed point in $X$.}\\

\noindent Lately, Zhang and Song [3, Theorem 2.1] (2009) proved a
theorem on the existence of a common fixed point for a duality of
two single valued generalized $\varphi$-weak contraction mappings.
By extending two  single valued mappings in the Theorem  of Zhang
and Song [3] to two multi-valued mappings, and By extending one
multi-valued mapping  in  Theorem 1.2 to a duality of multi-valued
mappings, Rouhani and Moradi [4] (2010) proved the following
coincidence theorem, without assuming $x\longrightarrow d(x, Tx)$ or
$x\rightarrow d(x, Sx)$ to be l.s.c.\\

\noindent{\bf Theorem 1.4} ([4, Theorem 3.1]). \textit{Let $(X, d)$
be a complete metric space, and let $T, S : X \rightarrow CB(X)$ be
two multivalued mappings such that for all $x, y\in X, H(Tx,Sy) \leq
\alpha M(x, y)$, where $0 \leq \alpha <1$. Then there exists a point
$x \in X$ such that $x \in Tx$ and $x \in Sx$ (i.e., $T$ and $S$
have a common fixed point). Moreover, if either $T$ or $S$ is single
valued, then
this common fixed point is unique.}\\\\
Further they also proved  the Theorem 1.5 below.\\

\noindent{\bf Theorem 1.5} ([4, Theorem 4.1]). \textit{Let $(X, d)$
be a complete metric space and let $T : X \rightarrow X$ and $S : X
\rightarrow CB(X)$ be two mappings such that for all $x, y \in X$,
$$H({Tx}, Sy)\leq M(x, y)-\varphi(M(x, y)),$$
 where $\varphi: [0,
+\infty)\rightarrow[0, +\infty)$ is l.s.c. with $\varphi(0)=0$ and
$\varphi(t)>0$ for all $t>0$. Then there exists a unique point
$x \in X$ such that $Tx = x \in Sx$.}\\

\noindent Finally, for the endpoint of multi-valued mappings,
Amini-Harandi [6] (2010) proved Theorem 1.6 below.\\

 \noindent{\bf
Theorem 1.6} ([6, Theorem 2.1]).\textit{Let $(X, d)$ be a complete
metric space and $T$ be a multi-valued mapping that satisfies
$$H(Tx, Ty)
\leq \varphi(d(x, y)),$$ for each $x, y \in X$ , where $\varphi :
[0,+\infty) \rightarrow [0,+\infty)$ is upper semicontinuous
(u.s.c.),  $\varphi(t) < t$ for each $t
> 0$  and satisfies $\liminf \limits_{
t\rightarrow \infty }(t - \varphi(t)) > 0$. Then $T$ has a unique
endpoint if and only if T has the approximate endpoint property.}\\

\noindent Moradi and Khojasteh [7] (2011) extended the result of
Amini-Harandi
to the following theorem 1.7.\\

\noindent{\bf Theorem 1.7} ([7, Theorem 2.1]). \textit{Let $(X, d)$
be a complete metric space and $T$ be a multi-valued mapping that
satisfies
$$H(Tx, Ty)
\leq \varphi(N(x, y)), $$ for each $  x, y \in X$ , where $\varphi :
[0,+\infty) \rightarrow [0,+\infty)$ is u.s.c. with $\varphi(t) < t$
for all $t
> 0$  and $\liminf \limits_{
t\rightarrow \infty }(t - \varphi(t)) > 0$. Then $T$ has a unique
endpoint if and only if T has the approximate endpoint property.}\\

Motivated by the contributions stated above, the present work make a
further study on the common fixed point for a duality of generalized
weak ($\varphi$-weak) contractions, and  also make a study on the
common endpoint for a duality of generalized weak ($\varphi$-weak)
contractions. Our contributions extend  the results of Theorem 1.3,
Theorem 1.4,  Theorem 1.5 and  Theorem 1.7.


\section{Preliminaries}
\label{1} This section proposes several Lemmas for our posterior
discussions.\\

\noindent{\bf Lemma 2.1.} \textit{Let $(X, d)$ be a complete metric
space
 and $S, T : X\rightarrow CB(X)$ are a duality of generalized weak
 (or
$\varphi$-weak) contractions. Then $Fix(S)=Fix(T)$.}\\

\noindent{\bf Proof}.  Let $x\in Fix(S)$. Then
$$\begin{array}{rcl}& &d(x, Tx)\leq H(Sx,Tx)\leq\alpha(x,x)M(x,x)\\
&=&\alpha(x,x)\max\{d(x,x), d(x,Sx), d(x,Tx), \frac{d(x,Tx)+d(x,Sx)
}{2}\}\\
&=&\alpha(x,x)d(x, Tx).\end{array}$$ Since $\alpha(x,x)<1$, this
implies $d(x,Tx)=0$.
That is, $x\in Fix(T)$. Hence, $Fix(S)=Fix(T)$.\\

\noindent{\bf Lemma 2.2.} \textit{Let $(X, d)$ be a complete metric
space, $\gamma\in[0,1)$,
and $x_{n}$ be a sequence of $X$  that satisfies
$$\begin{array}{rcl}d(x_n, x_{n+1}) \leq\gamma d(x_{n-1},
x_n)+\frac{1}{2^n}
\end{array}\eqno (2.1)$$
for all $n \in {\mathbb{N}}$ ($x_0 \in X$). Then  $\{x_n\}$ is
convergent.}\\

\noindent{\bf Proof}. By (2.1), for each $n\in {\mathbb{N}}$,
$$\begin{array}{rcl} d(x_n, x_{n+1})&\leq &\gamma d(x_{n-1},
x_n)+\frac{1}{2^n}\\&\leq&\gamma[\gamma d(x_{n-2},
x_{n-1})+\frac{1}{2^{n-1}}]+\frac{1}{2^n}\\&=&\gamma^2 d(x_{n-2},
x_{n-1})+\frac{\gamma}{2^{n-1}}+\frac{1}{2^n}\\& &\cdots\\&\leq&
\gamma^n d(x_{0}, x_{1})+\frac{\gamma^{n-1}}{2^{1}}+\cdots
+\frac{\gamma^1}{2^{n-1}}+\frac{\gamma^0}{2^n}\\&\leq&\frac{M}{1-\gamma}
(\frac{\gamma^n}{2^0}+\frac{\gamma^{n-1}}{2^{1}}+\cdots
+\frac{\gamma^1}{2^{n-1}}+\frac{\gamma^0}{2^n}),
\end{array}$$
where $M=\max\{d(x_{0}, x_{1}),1\}$. Without loss of generality,
assume $M=1$. Then
$$\begin{array}{rcl} d(x_n, x_{n+1})\leq
\frac{\gamma^n}{2^0}+\frac{\gamma^{n-1}}{2^{1}}+\cdots
+\frac{\gamma^1}{2^{n-1}}+\frac{\gamma^0}{2^n}.
\end{array}\eqno
(2.2)$$ By (2.2), for any $n, m \in {\mathbb{N}}$, we have
$$\begin{array}{rcl}& & d(x_n, x_{n+m})\\&\leq& d(x_n, x_{n+1})+d(x_{n+1},
x_{n+2})+\cdots+d(x_{n+m-1},
x_{n+m})\\&\leq&\{[\frac{\gamma^n}{2^0}+\frac{\gamma^{n-1}}{2^1}
+\frac{\gamma^{n-2}}{2^2}+\cdots+\frac{\gamma^{0}}{2^n}]\\&
&+[\frac{\gamma^{n+1}}{2^0}+\frac{\gamma^{n}}{2^1}
+\frac{\gamma^{n-1}}{2^2}+\cdots+\frac{\gamma^{1}}{2^n}+\frac{\gamma^{0}}{2^{n+1}}]\\&
&\cdots\\& &+ [\frac{\gamma^{n+m-1}}{2^0}+\frac{\gamma^{n+m-2}}{2^1}
+\frac{\gamma^{n+m-3}}{2^2}+\cdots\\&
&+\frac{\gamma^{m-1}}{2^n}+\frac{\gamma^{m-2}}{2^{n+1}}+\cdots+
\frac{\gamma^{0}}{2^{n+m-1}}]\}\\&=&\{[\frac{1}{2^0}(\gamma^{n}+\gamma^{n+1}+
\cdots+\gamma^{n+m-1})+\\& &\frac{1}{2^1}(\gamma^{n-1}+\gamma^{n}+
\cdots+\gamma^{n+m-2})+\cdots\\&
&+\frac{1}{2^n}(\gamma^{0}+\gamma^{1}+
\cdots+\gamma^{m-1})]+[\frac{1}{2^{n+1}}(\gamma^{0}+\gamma^{1}+
\cdots+\gamma^{m-2})\\& &+\frac{1}{2^{n+2}}(\gamma^{0}+\gamma^{1}+
\cdots+\gamma^{m-3})\\&
&+\cdots+\frac{1}{2^{n+m-2}}(\gamma^{0}+\gamma^{1})+\frac{1}{2^{n+m-1}}(\gamma^{0})]\}
\\&=&\{[\frac{1}{2^0}(\frac{\gamma^{n}-\gamma^{n+m}}{1-\gamma})+
\frac{1}{2^1}(\frac{\gamma^{n-1}-\gamma^{n+m-1}}{1-\gamma})+\cdots+
\frac{1}{2^n}(\frac{\gamma^{0}-\gamma^{m-1}}{1-\gamma})]\\&
&+[\frac{1}{2^{n+1}}(\frac{\gamma^{0}-\gamma^{m-1}}{1-\gamma})+
\frac{1}{2^{n+2}}(\frac{\gamma^{0}-\gamma^{m-2}}{1-\gamma})+\cdots+
\frac{1}{2^{n+m-1}}(\frac{\gamma^{0}-\gamma^{1}}{1-\gamma})]\}
\\&<&\frac{1}{(1-\gamma)}[(\frac{\gamma^{n}}{2^0}+
\frac{\gamma^{n-1}}{2^1}+\cdots+ \frac{\gamma^{0}}{2^n})\\&
&+(\frac{1}{2^{n+1}}+ \frac{1}{2^{n+2}}+\cdots+
\frac{1}{2^{n+m-1}})]\\&=
&\frac{1}{(1-\gamma)}\{\gamma^{n}[\frac{1}{(2\gamma)^0}+
\frac{1}{(2\gamma)^1}+\cdots+
\frac{1}{(2\gamma)^n}]+\frac{\frac{1}{2^{n+1}}-\frac{1}{2^{n+m}}}{1-\frac{1}{2}}\}
\\&<&\frac{1}{(1-\gamma)}\{\gamma^{n}[\frac{1-\frac{1}{(2\gamma)^{n+1}}}{1-\frac{1}{(2\gamma)}}]
+\frac{1}{2^{n}}\}.\end{array}\eqno (2.3)$$ In terms of (2.3), if
$(2\gamma)>1$, then
$$\begin{array}{rcl}& & d(x_n, x_{n+m})\\&<& \frac{1}{(1-\gamma)}
\{\gamma^n[\frac{2\gamma}{2\gamma-1}]+\frac{1}{2^n}\}=\frac{1}{(1-\gamma)}\cdot
\frac{(2\gamma)^{n+1}+2\gamma-1}{(2\gamma-1)2^n}\\&<&\frac{2\gamma
}{(2\gamma-1)(1-\gamma)}\cdot (\gamma^n+\frac{1}{2^n})<\frac{4
\gamma^{n+1}}{(2\gamma-1)(1-\gamma)}.\end{array}\eqno (2.4)$$
Otherwise, $(2\gamma)<1$, we have
$$\begin{array}{rcl}& & d(x_n, x_{n+m})\\&=& \frac{1}{(1-\gamma)}
[\gamma^n\cdot\frac{1-(2\gamma)^{n+1}}{(2\gamma)^{n}-(2\gamma)^{n+1}}+\frac{1}{2^n}]<
\frac{1}{(1-\gamma)}
[\frac{\gamma^n}{(2\gamma)^{n}-(2\gamma)^{n+1}}+\frac{1}{2^n}]\\&=&\frac{1}{(1-\gamma)}
[\frac{1}{2^{n}(1-2\gamma)}+\frac{1}{2^n}]=\frac{1}{(1-\gamma)(1-2\gamma)2^{n-1}}.\end{array}\eqno
(2.5)$$

From (2.4) and (2.5), we can easily know that the sequence $\{x_n\}$
is a Cauchy sequence. So it is convergent. This ends the proof.
$\square$\\

\noindent{\bf Lemma 2.3.} \textit{Let $(X, d)$ be a complete metric
space
 and $S, T : X\rightarrow CB(X)$ are a duality of generalized weak
contractions.  Let also $\{x_{n}\}$ be a convergent sequence of $X$
that satisfies $x_{n+1}\in Sx_{n}$ for each even $n \in
{\mathbb{N}}$, $\lim\limits_{n\rightarrow\infty}x_{n}=x^\ast$ and
$\limsup\limits_{k\rightarrow\infty}\alpha(x_{2k}, x^\ast)<1$. Then
$x^\ast\in Fix(T)=Fix(S)$.}\\

\noindent{\bf Proof}. In terms of the conditions, for each even $n
\in {\mathbb{N}}$, we have
$$\begin{array}{rcl}d(x_{n+1}, Tx^\ast) &\leq& H(Sx_{n}, Tx^\ast)
\leq\alpha(x_{n}, x^\ast)M(x_{n},
 x^\ast);
\end{array}\eqno
(2.6)$$
$$\begin{array}{rcl}& &M(x_{n}, x^\ast)\\&=&\max\{d(x_{n}, x^\ast), d(x_{n}, Sx_{n}),
d(x^\ast, Tx^\ast),\\& & \frac{d(x_{n}, Tx^\ast)+d(x^\ast,
Sx_{n})}{2}\}\\&\leq&\max\{d(x_{n}, x^\ast), d(x_{n}, x_{n+1}),
d(x^\ast, Tx^\ast),\\& & \frac{d(x_{n}, x^\ast)+d(x^\ast,
Tx^\ast)+d(x^\ast, x_{n})+d(x_{n},
Sx_{n})}{2}\}\\&\leq&\max\{d(x_{n}, x^\ast), d(x_{n}, x_{n+1}),
d(x^\ast, Tx^\ast),\\& &d(x_{n}, x^\ast)+ \frac{d(x^\ast,
Tx^\ast)+d(x_{n}, x_{n+1})}{2}\}.
\end{array}\eqno
(2.7)$$Note that $\lim\limits_{n\rightarrow\infty}x_{n}=x^\ast$.
Combining (2.6) and (2.7), we further obtain
$$\begin{array}{rcl}d(x^\ast, Tx^\ast)&\leq&[\limsup\limits_{k\rightarrow\infty}
\alpha(x_{2k}, x^\ast)]\limsup\limits_{k\rightarrow\infty}M(x_{2k},
x^\ast)
\\&\leq&[\limsup\limits_{k\rightarrow\infty}\alpha(x_{2k}, x^\ast)]d(x^\ast,
Tx^\ast).
\end{array}$$
 Since $\limsup\limits_{k\rightarrow\infty}\alpha(x_{2k}, x^\ast)<1$, this
implies $d(x^\ast, Tx^\ast)=0$. That is $x^\ast\in Tx^\ast$. The
proof
completes. $\square$\\

\noindent{\bf Lemma 2.4.} \textit{Let $(X, d)$ be a complete metric
space
 and $S, T : \rightarrow CB(X)$ are a duality of generalized weak
contractions. Then we have the the  conclusions as follows.\\
(1) $End(S)=End(T) (\subseteq Fix(S)=Fix(T))$ and $|End(S)|\leq 1$.
Here $|End(S)|$ denotes the cardinal number of $End(S)$. (This
implies that $S$ and $T$ have an unique common endpoint, or  have no
endpoint.)\\
(2) If $S$ and  $T$ have common endpoint, then $\inf\limits_{x\in
X}[H(\{x\}, Sx)+H(\{x\}, Tx)]$=0, termed as   the
approximate endpoint property of  duality $S$ and $T$.\\
(3) If either $S$ or  $T$ is single valued, then
$End(S)=End(T)=Fix(S)=Fix(T)$. (This implies that the fixed points
of $S$ and
$T$ must be endpoints.)}\\

\noindent{\bf Proof}.   Let $x\in End(S)$. Then  $x\in
Fix(S)=Fix(T)$ from Lemma 2.1. This implies $M(x,x)=0$. Therefore,
we have
$$\begin{array}{rcl}& &H(\{x\}, Tx)= H(Sx,Tx)\leq\alpha(x,x)M(x,x)=0.\end{array}$$
This means $Tx=\{x\}$. That is, $x\in End(T)$. Hence
$End(S)=End(T)$.

Let $x,y\in End(S)=End(T)$. Then $M(x,y)=d(x,y)$, further
$$d(x,y)=H(\{x\}, \{y\})=
H(Sx,Ty)\leq\alpha(x,y)M(x,y)=\alpha(x,y)d(x,y).$$ For
$\alpha(x,y)<1$, this implies $d(x,y)=0$. That is $x=y$. Hence
$|End(S)|\leq 1$.

We have proved (1). (2) is obvious. Next we further prove (3).

Suppose that one of $S$ and $T$ is single valued. Without loss of
generality, we assume $S$ is single valued. Then it is obvious that
$End(S)=Fix(S)$. So $End(T)=End(S)=Fix(S)=Fix(T)$. This ends the
proof. $\square$

\section{Fixed
point theory} \label{1} In the  section, we focus on studying the
fixed point theory.

We are now in a position to prove our first theorem, which extends
Theorem 2.3 of  Daffer and Kaneko [1] by generalizing one mapping
$T$ to two  mappings $S$ and $T$, and by improving the other
conditions, which also extends Theorem 3.1 of Rouhani and Moradi [4]
by replacing the constant contraction factor $\alpha$ with an
general
$\alpha(x, y)$.\\

\noindent{\bf Theorem 3.1.} \textit{Let $(X, d)$ be a complete
metric space
 and $S, T : \rightarrow CB(X)$ are a duality of generalized weak
contractions that satisfies
$$\sup\{\alpha(x_{2k-2},
 x_{2k-1}),\alpha(x_{2k},
 x_{2k-1})|k\in{\mathbb{N}}\}<1 \eqno (3.1)$$
for any  sequence  $\{x_n\}$ of $X$ with  $\{d(x_{n}, x_{n+1})\}$ to
be monotone decreasing, and $\alpha$ is
 u.s.c. (or
$\limsup\limits_{n\rightarrow\infty}\alpha(x_{n}, x^\ast)<1$ if
$\lim\limits_{n\rightarrow\infty}x_{n}=x^\ast$).
Then $Fix(S)=Fix(T)\neq\emptyset$.}\\

\noindent{\bf Proof.} (1) By Lemma 2.1, $Fix(S)=Fix(T)$. To complete
the proof, what we need is only to prove
$Fix(S)=Fix(T)\neq\emptyset$. Arguing by contradiction, we assume
$Fix(S)=Fix(T)=\emptyset$.

\noindent(2) Let $x_0 \in X$. Then $d(x_0,Sx_0)>0$. It is obvious
that we can choose a $x_1 \in Sx_0$ such that $0<d(x_0,
x_1)<d(x_0,Sx_0)+1$, and $d(x_1,Tx_1)>0$.

Let $\varepsilon_1=\min\{\frac{1}{2}, [1-\alpha(x_0, x_1)]d(x_0,
x_1)\}$. Then there exists a $x_2 \in Tx_1$ such that $0<d(x_1,
x_2)<d(x_1,Tx_1)+\varepsilon_1$. Let
$\varepsilon_2=\min\{\frac{1}{2^2}, [1-\alpha(x_2, x_1)]d(x_1,
x_2)\}$. Then there exists a $x_3 \in Sx_2$ such that $0<d(x_2,
x_3)<d(x_2,Sx_2)+\varepsilon_2$. Inductively, we have the  general
fact as follows.
 For each $k\in {\mathbb{N}}$, let
$$\varepsilon_{2k-1}=\min\{\frac{1}{2^{2k-1}}, [1-\alpha(x_{2k-2},
x_{2k-1})]d(x_{2k-2}, x_{2k-1})\}.$$ Then there exists a $x_{2k} \in
Tx_{2k-1}$ such that
$$0<d(x_{2k-1},
x_{2k})<d(x_{2k-1},Tx_{2k-1})+\varepsilon_{2k-1}\leq
d(x_{2k-1},Tx_{2k-1})+\frac{1}{2^{2k-1}}.$$ Let also
$$\varepsilon_{2k}=\min\{\frac{1}{2^{2k}}, [1-\alpha(x_{2k},
x_{2k-1})]d(x_{2k-1}, x_{2k})\}.$$ Then there exists a $x_{2k+1} \in
Sx_{2k}$ such that $$0<d(x_{2k},
x_{2k+1})<d(x_{2k},Sx_{2k})+\varepsilon_{2k}\leq
d(x_{2k},Sx_{2k})+\frac{1}{2^{2k}}.$$

\noindent(3) For the sequence $\{x_{n}\}$ constructed above,
$\forall n \in {\mathbb{N}}$, when $n$ is odd, we have
$$[1-\alpha(x_{n-1}, x_{n})]d(x_{n-1}, x_{n})\geq
\varepsilon_{n}.\eqno (3.2)$$ Further,
$$\begin{array}{rcl}d(x_n, x_{n+1})&< &
d(x_n, Tx_{n})+\varepsilon_{n} \leq H(Sx_{n-1},
Tx_{n})+\varepsilon_{n} \\&\leq&\alpha(x_{n-1}, x_n)M(x_{n-1},
x_n)+\varepsilon_{n}\\&\leq&\alpha(x_{n-1}, x_n)M(x_{n-1},
x_n)+\frac{1}{2^n},
\end{array}\eqno (3.3)$$
$$\begin{array}{rcl} & &
 M(x_{n-1}, x_n)\\
& \leq& \max\{d(x_{n-1}, x_n), d(x_{n-1}, Sx_{n-1}), d(x_n,
Tx_n),\\& & \frac{d(x_{n-1}, Tx_n) + d(x_{n}, Sx_{n-1})}{2}\}\\&\leq
& \max\{d(x_{n-1}, x_n), d(x_{n-1}, x_{n}), d(x_n, x_{n+1}),\\& &
\frac{d(x_{n-1}, x_n) +d(x_{n}, Tx_{n})}{2}\} \\& \leq&
\max\{d(x_{n-1}, x_n), d(x_n, x_{n+1}), \frac{d(x_{n-1}, x_n)
+d(x_{n}, x_{n+1})}{2}\}\\  & =& \max\{d(x_{n-1}, x_n), d(x_n,
x_{n+1})\}.
\end{array}\eqno (3.4)$$
If $\max\{d(x_{n-1}, x_n), d(x_n, x_{n+1})\}=d(x_n, x_{n+1})$, i.e.
$d(x_{n-1}, x_n)\leq d(x_n, x_{n+1})$, then
$$[1- \alpha(x_{n-1}, x_n)]d(x_{n-1}, x_n)\leq[1- \alpha(x_{n-1},
x_n)]d(x_{n}, x_{n+1}),\eqno (3.5)$$ and  from (3.3) and (3.4), we
obtain
$$\begin{array}{rcl}& &d(x_n, x_{n+1})<
\alpha(x_{n-1}, x_n)d(x_n, x_{n+1})+\varepsilon_{n}\\&\Rightarrow&
[1- \alpha(x_{n-1}, x_n)]d(x_{n}, x_{n+1})<\varepsilon_{n}.
\end{array}\eqno (3.6)$$
Combing (3.5) and (3.6), we obtain  $[1- \alpha(x_{n-1},
x_n)]d(x_{n-1}, x_n)<\varepsilon_{n}$. This contradicts (3.2).So
$\max\{d(x_{n-1}, x_n), d(x_n, x_{n+1})\}\neq d(x_{n-1}, x_n)$.This
yields to $\max\{d(x_{n-1}, x_n), d(x_n, x_{n+1})\}=d(x_{n-1}, x_n)$
and $d(x_n, x_{n+1})<d(x_{n-1}, x_n)$. Also, from (3.3) and (3.4),
we obtain
$$\begin{array}{rcl}d(x_n, x_{n+1})&< &
\alpha(x_{n-1}, x_n)d(x_{n-1}, x_n)+\frac{1}{2^n}.
\end{array}\eqno (3.7)$$
When $n$ is even, we have
$$[1-\alpha(x_{n}, x_{n-1})]d(x_{n-1}, x_{n})\geq
\varepsilon_{n},\eqno (3.8)$$
$$\begin{array}{rcl}& &d(x_n, x_{n+1})=d( x_{n+1},x_n)\\&< &
d(Sx_n, x_{n})+\varepsilon_{n} \leq H(Sx_{n},
Tx_{n-1})+\varepsilon_{n} \\&\leq&\alpha( x_n,x_{n-1})M(x_{n},
x_{n-1})+\varepsilon_{n}\\&\leq&\alpha( x_n,x_{n-1})M(x_{n},
x_{n-1})+\frac{1}{2^n},
\end{array}\eqno (3.9)$$
$$\begin{array}{rcl} & &
M(x_{n}, x_{n-1})\\
& \leq& \max\{d(x_{n-1}, x_n), d(x_{n}, Sx_{n}), d(x_{n-1},
Tx_{n-1}),\\& & \frac{d(x_{n}, Tx_{n-1}) + d(x_{n-1},
Sx_{n})}{2}\}\\&\leq & \max\{d(x_{n-1}, x_n), d(x_n, x_{n+1}),
d(x_{n-1}, x_{n}),\\& & \frac{d(x_{n-1}, x_n) +d(x_{n},
Sx_{n})}{2}\}
\\& \leq& \max\{d(x_{n-1}, x_n), d(x_n, x_{n+1}), \frac{d(x_{n-1},
x_n) +d(x_{n}, x_{n+1})}{2}\}\\  & =& \max\{d(x_{n-1}, x_n), d(x_n,
x_{n+1})\}.
\end{array}\eqno (3.10)$$
From $(3.8)$, (3.9) and (3.10), in the same way as used above, we
can also obtain $d(x_n, x_{n+1})<d(x_{n-1}, x_n)$ and
$$\begin{array}{rcl}d(x_n, x_{n+1})&< &
\alpha(x_n, x_{n-1})d(x_{n-1}, x_n)+\frac{1}{2^n}.
\end{array}\eqno (3.11)$$
(4)
From (3), it is obvious that the sequence $\{d(x_{n}, x_{n+1})\}$ is
monotone decreasing. Hence  (3.1) holds. So, there exists a
$\gamma<1$ such that $$\max\{\alpha(x_{2k-2},
 x_{2k-1}),\alpha(x_{2k},
 x_{2k-1})\}<\gamma$$ for all $k\in {\mathbb{N}}$. Therefore
using (3.7) and $(3.11)$ we can obtain (2.1). Thus $\{x_{n}\}$ is
convergent from Lemma 2.2.

Finally, let  $\lim\limits_{n\rightarrow\infty}x_{n}=x^\ast$. Then,
since $\alpha$ is u.s.c. we have
$\limsup\limits_{n\rightarrow\infty}\alpha(x_{n},
x^\ast)\leq\alpha(x^\ast, x^\ast)<1$. Note that the approach we
produce the sequence $\{x_{n}\}$. By Lemma 2.3, $x^\ast\in Tx^\ast$.
This contradicts $Fix(T)=\emptyset$. So
$Fix(S)=Fix(T)\neq\emptyset$. The proof completes. $\square$

 As an application we  propose a proof of Theorem 1.3 from Theorem
 3.1
 as follows.

\noindent{\bf Proof of Theorem 1.3}. let
$$\begin{array}{rcl}\alpha(x, y)=
\left\{\begin{array}{l}1-\frac{\phi(x, y)}{d(x, y)}, d(x, y)\neq
0;\\0, d(x, y)=0,\end{array}\right.
\end{array}$$
and $S=T$. Then $S$ and $T$ are a duality of generalized weak
contractions. Let also $\{x_{n}\}$ be a  sequence of $X$ with
$\{d(x_{n},x_{n+1})\}$ to be monotone decreasing. And assume
$\lim\limits_{n\rightarrow \infty}d(x_{n},x_{n+1})=r$.

If $r>0$, then $\lambda(r,d(x_{1},x_{2}))>0$ for $\phi$ is compactly
positive. On the other hand,
$\phi(x_{n},x_{n+1})\geq\lambda(r,d(x_{1},x_{2}))$ since
$r<d(x_{n},x_{n+1})\leq d(x_{1},x_{2})$. So we have
$$\alpha(x_{n},x_{n+1})=1-\frac{\phi(x_{n},x_{n+1})}{d(x_{n},x_{n+1})}\leq
1-\frac{\lambda(r,d(x_{1},x_{2}))}{d(x_{1},x_{2})}<1.$$ With the
same argument, $\alpha(x_{n+1},x_{n})\leq
1-\frac{\lambda(r,d(x_{1},x_{2}))}{d(x_{1},x_{2})}$. Hence (3.1)
holds. If $r=0$, then
$$\begin{array}{rcl}\limsup\limits_{n\rightarrow
\infty}\alpha(x_{n},x_{n+1})&=&\limsup\limits_{n\rightarrow
\infty}[1-\frac{\phi(x_{n},x_{n+1})}{d(x_{n},x_{n+1})}]\\&\leq&
\limsup\limits_{n\rightarrow
\infty}[1-\frac{\lambda(d(x_{n-1},x_{n}),
d(x_{n},x_{n+1}))}{d(x_{n},x_{n+1})}]\\&\leq&
1-\liminf\limits_{n\rightarrow
\infty}\frac{\lambda(d(x_{n-1},x_{n}),
d(x_{n},x_{n+1}))}{d(x_{n},x_{n+1})}\\&\leq&1-\liminf\limits_{\beta\rightarrow
\infty}\frac{\lambda(\alpha, \beta)}{\beta}<1.\end{array}$$ With the
same argument, $\limsup\limits_{n\rightarrow
\infty}\alpha(x_{n+1},x_{n})<1$. Hence (3.1) holds.

Let $\lim\limits_{n\rightarrow \infty}x_{n}=x^\ast$. Without loss of
generality, assume $d(x_{n}, x^\ast)\neq 0$ for all $
n\in{\mathbb{N}}$. Then
$$\begin{array}{rcl}& &\limsup\limits_{n\rightarrow\infty}\alpha(x_{n}, x^\ast)
=\limsup\limits_{n\rightarrow\infty}[1-\frac{\phi(x_{n},
x^\ast)}{d(x_{n}, x^\ast)}]\leq\\& &
1-\liminf\limits_{n\rightarrow\infty}\frac{\phi(x_{n},
x^\ast)}{d(x_{n}, x^\ast)}\leq 1-\liminf\limits_{\beta\rightarrow
0}\frac{\lambda(\alpha, \beta)}{\beta}<1.\end{array}$$

Combing the results above, by Theorem 3.1,  $T$ has fixed point.
This ends the proof.
$\square$\\

\noindent{\bf Theorem 3.2.} \textit{Let $(X, d)$ be a complete
metric space
 and $S, T : \rightarrow CB(X)$ are a duality of generalized $\varphi$-weak
contractions that satisfies $\varphi$ is u.s.c. and
$$\limsup\limits_{t\rightarrow 0}\frac{\varphi(t)}{t}<1.\eqno (3.13)$$
Then $Fix(S)=Fix(T)\neq\emptyset$.}\\

\noindent{\bf Proof}. For any  $(x, y)\in X\times X$, put
$$\begin{array}{rcl}\alpha(x, y)=
\left\{\begin{array}{l}\frac{\varphi(M(x, y))}{M(x, y)}, M(x, y)\neq
0;\\0, M(x, y)=0.\end{array}\right.
\end{array}$$
Then it can be easily verify that $H(Sx, Ty)\leq\alpha(x, y)M(x,
y)$. That is, $S, T : \rightarrow CB(X)$ are  a duality of
generalized weak contractions with the $\alpha(x, y)$. Note that the
conditions $\alpha$ is u.s.c. and (3.1) are used only in the step
(4) of the proof of Theorem 3.1. We can easily know that the steps
(1), (2) and (3) can be used to prove Theorem 3.2. So the proof can
be accomplished by proposing the step $(4)'$ below.

\noindent $(4)'$ $\forall n\in{\mathbb{N}}$, assume first that $n$
is odd. Note that $0<d(x_{n-1}, x_{n})\leq M(x_{n-1}, x_{n})$ and
$$\max\{d(x_{n-1}, x_{n}), d(x_n, x_{n+1})\}=d(x_{n-1}, x_{n}).$$ By
(3.4), we have $M(x_{n-1}, x_{n})=d(x_{n-1}, x_{n})>0$. This leads
to
$$\alpha(x_{n-1}, x_{n})=\frac{\varphi(d(x_{n-1}, x_{n}))}{d(x_{n-1},
x_{n})}.\eqno (3.14)$$ Further, from (3.7), we obtain $$d(x_{n},
x_{n+1})\leq \varphi(d(x_{n-1}, x_{n}))+\frac{1}{2^n}.\eqno (3.15)$$
When  $n$ is even, with the same argument, we have (3.15) and
$$\alpha(x_{n}, x_{n-1})=\frac{\varphi(d(x_{n-1},
x_{n}))}{d(x_{n-1}, x_{n})}.\eqno (3.16).$$ Since the sequence
$\{d(x_{n}, x_{n+1})\}$ is monotone decreasing and bounded below, it
is convergent. Let $\lim\limits_{n\rightarrow\infty}d(x_n,
x_{n+1})=r$. For $\varphi$ is u.s.c. using (3.15) we have $r\leq
\varphi(r)$. This implies $r=0$ because $\varphi(t)<t$ for all
$t>0$. Therefore, according to (3.14) and (3.16), we respectively
have
$$\limsup\limits_{k\rightarrow\infty}\alpha(x_{2k-2}, x_{2k-1})=
\limsup\limits_{k\rightarrow\infty}\frac{\varphi(x_{2k-2},
x_{2k-1})}{d(x_{2k-2}, x_{2k-1})}\leq \limsup\limits_{t\rightarrow
0}\frac{\varphi(t)}{t}<1.$$
$$\limsup\limits_{k\rightarrow\infty}\alpha(x_{2k}, x_{2k-1})=
\limsup\limits_{k\rightarrow\infty}\frac{\varphi(x_{2k},
x_{2k-1})}{d(x_{2k}, x_{2k-1})}\leq \limsup\limits_{t\rightarrow
0}\frac{\varphi(t)}{t}<1.$$Hence (3.1) holds.
 Using
 (3.7) and $(3.11)$ we obtain (2.1). Thus $x_{n}$ is convergent from Lemma 2.2.

Finally, let  $\lim\limits_{n\rightarrow\infty}x_{n}=x^\ast$. Then
for each even $n$, we have (2.7). This reduces to
$\limsup\limits_{k\rightarrow\infty}M(x_{2k}, x^\ast)\leq d(x^\ast,
Tx^\ast)$. So there exists a positive number $b$ such that
$M(x_{2k}, x^\ast)\leq b$. For $\varphi$ is u.s.c. and
$\limsup\limits_{t\rightarrow 0}\frac{\varphi(t)}{t}<1$,
$\sup\{\frac{\varphi(t)}{t}|t\in(0, b]\}<1$. Therefore,
$\limsup\limits_{n\rightarrow\infty}\alpha(x_{n},
x^\ast)=\limsup\limits_{n\rightarrow\infty}\frac{\varphi(M(x_{n},
x^\ast))}{M(x_{n}, x^\ast)}\leq\sup\{\frac{\varphi(t)}{t}|t\in(0,
b]\}<1$. By Lemma 2.3, $x^\ast\in Tx^\ast$. This contradicts
$Fix(T)=\emptyset$. So $Fix(S)=Fix(T)\neq\emptyset$. The proof ends.
$\square$

Theorem 3.2 extends Theorem 4.1 of Rouhani and Moradi [4] by
allowing both  two mappings $S$ and $T$ to be multi-valued. However,
we add the condition (3.1). Whether Theorem 3.2 holds or not without
the condition (3.1) is a topic for us to further pursue.

\section{Endpoint theory} \label{1}Now we turn to address the endpoint theory.

In terms of Theorem 3.1 (Theorem 3.2) and Lemma 2.4, we can
immediately get the next corollary.\\

{\bf corollary 3.2.} \textit{Under the conditions of Theorem 3.1 (
Theorem 3.2), if  either $S$ or  $T$ is single valued, then there
exists a unique common fixed  point for $S$ and  $T$, which is also
a unique common endpoint of theirs.}\\

For $S$ and $T$ are all  multi-valued, we have Theorem 4.1 and
Theorem 4.2 below.\\

{\bf Theorem 4.1.} \textit{Let $(X, d)$ be a complete metric space
 and $S, T : X\rightarrow CB(X)$ are  a duality of generalized weak
contractions that satisfies $\alpha$ is u.s.c. and
$$\limsup\limits_{n, m\rightarrow\infty}\alpha(x_n,x_m)<1\eqno (4.1)$$
if $\lim\limits_{n,
m\rightarrow\infty}d(x_n,x_m)[1-\alpha(x_n,x_m)]=0$. Then $S$ and
$T$ have a unique common endpoint if they have the approximate
endpoint property.}\\

\noindent{\bf Proof}.  Suppose that $S$ and $T$ have the approximate
endpoint property. Then there exists a sequence $\{x_n\}$ such that
$$\lim\limits_{n\rightarrow \infty}[H(\{x_n\}, Sx_n)+ H(\{x_n\},
Tx_n)]= 0.$$

For all $m, n \in {\mathbb{N}}$, we have
$$\begin{array}{rcl}& & M(x_n, x_m)\\ &=& \max\{d(x_n, x_m), d(x_n, Sx_n), d(x_m,
Tx_m), \\& &\frac{d(x_n, Tx_m) + d(x_m, Sx_n)}{ 2 }\}\\& \leq &
\max\{d(x_n, x_m), H(\{x_n\}, Sx_n), H(\{x_m\}, Tx_m),\\& &
\frac{d(x_n, x_m)+H(\{x_m\}, Tx_m) +d(x_n, x_m)+ H(\{x_n\},
Sx_n)}{2}\} \\&\leq & d(x_n, x_m) + H(\{x_n\}, Sx_n) + H(\{x_m\},
Tx_m). \end{array}\eqno (4.2)$$ Note that $d(x_n, x_m)\leq
H(\{x_n\}, Sx_n) +H(Sx_n, Tx_m)+ H(\{x_m\}, Tx_m)$. From (4.2), we
further have
$$\begin{array}{rcl}& & M(x_n, x_m)\\&\leq &d(x_n, x_m)
-H(\{x_n\}, Sx_n) - H(\{x_m\}, Tx_m)\\& & + 2H(\{x_n\}, Sx_n) +
2H(\{x_m\}, Tx_m)\\ &\leq & H(Tx_n, Tx_m) + 2H(\{x_n\}, Sx_n) +
2H(\{x_m\}, Tx_m). \end{array}\eqno (4.3)$$ This reduces to
$$\begin{array}{rcl}& & M(x_n, x_m)\\&\leq & \alpha(x_n,
x_m) M(x_n, x_m) + 2H(\{x_n\}, Sx_n) + 2H(\{x_m\},
Tx_m).\end{array}\eqno (4.4)$$ Note that $d(x_n, x_m)\leq  M(x_n,
x_m)$. Using (4.4) we obtain
$$\begin{array}{rcl}& & d(x_n, x_m)[1-\alpha(x_n,
x_m)]\\&\leq & M(x_n, x_m)[1-\alpha(x_n, x_m)]\\&\leq & 2H(\{x_n\},
Sx_n) + 2H(\{x_m\}, Tx_m)\\&\Rightarrow &\lim\limits_{
n,m\rightarrow \infty}d(x_n, x_m)[1-\alpha(x_n, x_m)]=
0.\end{array}\eqno (4.5)$$ For (4.5), we have (4.1). Using also
(4.4), we obtain
$$\limsup\limits_{
n,m\rightarrow \infty}M(x_n, x_m)\leq[\limsup\limits_{
n,m\rightarrow \infty}\alpha(x_n, x_m)]\limsup\limits_{
n,m\rightarrow \infty}M(x_n, x_m).$$ By (4.1), this yields to
$\limsup\limits_{ n,m\rightarrow \infty}M(x_n, x_m)=0$. Hence
$\limsup\limits_{ n,m\rightarrow \infty}d(x_n, x_m)=0$, i.e.
$\{x_n\}$ is a Cauchy sequence.

Let $\lim\limits_{ n\rightarrow \infty}x_n=x^\ast$. For all $ n \in
{\mathbb{N}}$, we have
$$\begin{array}{rcl}& &H(\{x_n\}, Tx^\ast)- H(\{x_n\}, Sx_n)\\&\leq& H(Sx_n, Tx^\ast)
\leq \alpha(x_n, x^\ast)M(x_n, x^\ast)\\& =&\alpha(x_n,
x^\ast)\max\{d(x_n, x^\ast), d(x_n, Sx_n), d(x^\ast, Tx^\ast),\\&
&\frac{d(x_n, Tx^\ast)+d(x^\ast, Sx_n)}{2}\}\\ &\leq&\alpha(x_n,
x^\ast)\max\{d(x_n, x^\ast), d(x_n, Sx_n), d(x^\ast, Tx^\ast),\\&
&\frac{d(x_n, x^\ast)+d(x^\ast, Tx^\ast)+d(x^\ast, x_n)+d(x_n,
Sx_n)}{2}\}.
\end{array}$$Noting also $\alpha$ is u.s.c. we obtain
$$H(\{x^\ast\}, Tx^\ast)\leq\alpha(x^\ast,
x^\ast)d(x^\ast, Tx^\ast)\leq\alpha(x^\ast, x^\ast)H(\{x^\ast\},
Tx^\ast).$$ For $\alpha(x^\ast, x^\ast)<1 $, we conclude that
$H({x^\ast}, Tx^\ast)=0$. This means  $Tx^\ast = \{x^\ast\}$.
Finally, the uniqueness of the endpoint is concluded from  Lemma
2.4. $\square$

The following Theorem 4.2 is our final result, which extends the
Theorem 2.1 of Moradi and Khojasteh [7] to the case where both two
mappings are  multi-valued.\\

{\bf Theorem 4.2.}  \textit{Let $(X, d)$ be a complete metric space
 and $S, T : \rightarrow CB(X)$ are a duality of $\varphi$-generalized weak
contractions that satisfies $\varphi$ is s.u.c. and
$$\liminf\limits_{t\rightarrow \infty}[t-\varphi(t)]>0.\eqno (4.6)$$
Then $S$ and  $T$ have a unique common endpoint if they have the
approximate endpoint property.}\\

\noindent{\bf Proof}. Suppose that $S$ and $T$ have the approximate
endpoint property. Then there exists a sequence $\{x_n\}$ such that
$$\lim\limits_{n\rightarrow \infty}[H(\{x_n\}, Sx_n)+ H(\{x_n\},
Tx_n)]= 0,$$as well as (4.2) and (4.3) hold.

By (4.3) we  have
$$\begin{array}{rcl}& & M(x_n, x_m)\\&\leq & \varphi( M(x_n, x_m)) + 2H(\{x_n\}, Sx_n) + 2H(\{x_m\},
Tx_m).\end{array}\eqno (4.7)$$ If $\limsup\limits_{ n,m\rightarrow
\infty}M(x_n, x_m)=+\infty$, then
$$\begin{array}{rcl}\liminf\limits_{
t\rightarrow \infty}[t-\varphi(t)]\leq\liminf\limits_{
n,m\rightarrow \infty}[M(x_n, x_m)-\varphi(M(x_n, x_m))]\leq
0.\end{array}$$ This contradicts (4.6). So $\limsup\limits_{
n,m\rightarrow \infty}M(x_n, x_m)<+\infty$. Noting also $\varphi(t)$
is u.s.c. and using (4.7), we obtain
$$\limsup\limits_{
n,m\rightarrow \infty}M(x_n, x_m)\leq\limsup\limits_{ n,m\rightarrow
\infty}\varphi(M(x_n, x_m))\leq\varphi(\limsup\limits_{
n,m\rightarrow \infty}M(x_n, x_m)).$$ Note that $\limsup\limits_{
n,m\rightarrow \infty}M(x_n, x_m)<+\infty$ and $\varphi(t)<t$ for
all $t>0$. This implies
 $\limsup\limits_{ n,m\rightarrow
\infty}M(x_n, x_m)=0$. Thus  $\{x_n\}$ is Cauchy sequence.

Let $\lim\limits_{ n\rightarrow \infty}x_n=x^\ast$. For all $ n \in
{\mathbb{N}}$, we have
$$\begin{array}{rcl}& &H(\{x_n\}, Tx^\ast)- H(\{x_n\}, Sx_n)\\&\leq& H(Sx_n, Tx^\ast)
\leq \varphi(M(x_n, x^\ast)).
\end{array}$$ This reduces to
$$H(\{x^\ast\}, Tx^\ast)\leq \limsup\limits_{ n,m\rightarrow \infty}\varphi(M(x_n,
x^\ast))\leq\varphi(\limsup\limits_{ n,m\rightarrow \infty}M(x_n,
x^\ast)). \eqno(4.8)$$On the other hand,
$$\begin{array}{rcl}M(x_n, x^\ast) &\leq&\max\{d(x_n,
x^\ast), d(x_n, Sx_n), d(x^\ast, Tx^\ast),\\& &\frac{d(x_n,
x^\ast)+d(x^\ast, Tx^\ast)+d(x^\ast, x_n)+d(x_n, Sx_n)}{2}\}.
\end{array}\eqno(4.9)$$
If $H(\{x^\ast\}, Tx^\ast)\neq 0$, from (4.8) and (4.9), we have
$$\begin{array}{rcl}& &H(\{x^\ast\}, Tx^\ast)<\limsup\limits_{ n,m\rightarrow \infty}M(x_n,
x_m)\leq d(x^\ast, Tx^\ast)\leq H(\{x^\ast\}, Tx^\ast).\end{array}$$
This contradiction shows $H(\{x^\ast\}, Tx^\ast)=0$.  That is,
$Tx^\ast = \{x^\ast\}$. Finally, the uniqueness of the endpoint is
concluded from  Lemma 2.4. $\square$

\noindent{\bf Remark 4.3.}  By taking $S=T$, we
 can immediately obtain the Theorem 2.1 of Moradi and Khojasteh [7]
  from Theorem 4.2 and Lemma
2.4.\\

{\bf Acknowledgements}

The author cordially thank the anonymous referees for their valuable
comments which lead to the improvement of this paper.





\bibliographystyle{model1-num-names}
\bibliography{<your-bib-database>}







\end{document}